\documentclass{article}
\usepackage[english]{babel}
\usepackage[a4paper,top=2cm,bottom=2cm,left=3cm,right=3cm,marginparwidth=1.75cm]{geometry}
\usepackage{amsmath, amsthm, amsfonts, amssymb, amscd}
\usepackage{mathrsfs}
\usepackage{graphicx}
\usepackage{color}
\usepackage[colorlinks=true,hypertexnames=false,linkcolor=blue,citecolor=blue, backref=page]{hyperref}
\usepackage{tikz,pgfplots}
\usetikzlibrary{arrows}
\usepackage{newclude}
\usepackage{setspace}
\usepackage{enumitem}
\usepackage{latexsym}
\usepackage{verbatim}
\usepackage{cite}
\usepackage[noblocks]{authblk}
\usepackage[utf8]{inputenc}
\pdfoutput=1
\advance \textheight by \topskip
\numberwithin{equation}{section}
\newtheorem{teo}{Theorem}
\newtheorem*{teo*}{Theorem}

\usepackage{fancyhdr}
\makeatletter
\renewenvironment{proof}[1][\proofname]{\par
\pushQED{\qed}%
\normalfont \topsep6\p@\@plus6\p@\relax
\trivlist
\item\relax
{\bfseries
#1\@addpunct{.}}\hspace\labelsep\ignorespaces
}{%
\popQED\endtrivlist\@endpefalse
}
\makeatother
\def\ackname{Acknowledgements}%
\newenvironment{ack}[1][\ackname]%
{\ifx#1\empty\else\subsection*{#1.}\fi\par}
{\par}
\usepackage{cancel}
\usepackage[normalem]{ulem}
\newcommand\redst{\bgroup\markoverwith{\textcolor{red}{\rule[0.5ex]{2pt}{0.7pt}}}\ULon}
\newcommand\bluest{\bgroup\markoverwith{\textcolor{blue}{\rule[0.5ex]{2pt}{0.7pt}}}\ULon}

\usepackage[textwidth=2.7cm, color=red!20]{todonotes}
\newcounter{nota}[section]

\usepackage{xparse}
\NewDocumentCommand{\hess}{om}{%
\IfNoValueTF{#1}
{{\rm Hess} (#2)} {{\rm Hess} (#2) \Big\lvert _{#1}} 
}
\usepackage{yhmath}
\usepackage{mathtools}

\newcommand\restric[2]{ {
	\left.\kern-\nulldelimiterspace 
	#1 
	\vphantom{\big|} 
	\right|_{#2}
	}}
\makeatletter
\let\@fnsymbol\@arabic
\makeatother

\newcommand{\dr}{{{\tt x}_r}}
\newcommand{\dtheta}{{{\tt x}_{\theta}}}

\title{Optical Billiards}
\date{\today \\ 
\small 
Celebrating the 80th birthday of Marco Ant\^onio Teixeira
}
\author{S. Oliffson Kamphorst\thanks{Departamento de Matem\'atica, UFMG, Brasil},
S. Pinto-de-Carvalho$^1$
and C. H. Vieira Morais\thanks{Departamento de Matem\'atica, UFES, Brasil}
}
\begin{document}
\maketitle

\begin{abstract}
This paper investigates the dynamics of optical billiards, a generalization of classic billiards, where light rays travel within a refractive medium and reflect elastically at the boundary \cite{astro-spc}.
Inspired by studies of acoustic modes in rapidly rotating stars \cite{lig},
optical billiards are approached using Riemannian geometry, treating the refractive index as a metric.
The discrete dynamics is explored through geodesics, curvature properties, highlighting similarities and 
differences with classical billiards.
Examples include geodesic circular billiards and refractive media with zero Gaussian curvature,
revealing integrability conditions and singular behaviors.
The study establishes foundational results for optical billiards, emphasizing their dependence on boundary convexity and refractive index. 
\end{abstract}


\section{Introduction}

In \cite{lig}, Ligni\`eres and Georget present an asymptotic analysis of high-frequency acoustic modes in rapidly rotating stars, based on acoustic 
ray dynamics. They observe that, as the star rotates, its originally spherical boundary deforms in the equatorial plane, transforming into a symmetric ovoid. Using rotational symmetries, they show that the acoustic rays can be described 
using two ingredients:
a classic Hamiltonian of two degrees of freedom depending on the star's rotation frequency with a potential tending to infinity as it approaches the star's boundary and an isotropic two-dimensional medium, with refractive index depending on the distance to the center of the star.

\begin{figure}[ht]
\begin{center}
 \parbox{.14\hsize}{
\includegraphics[width=\hsize]{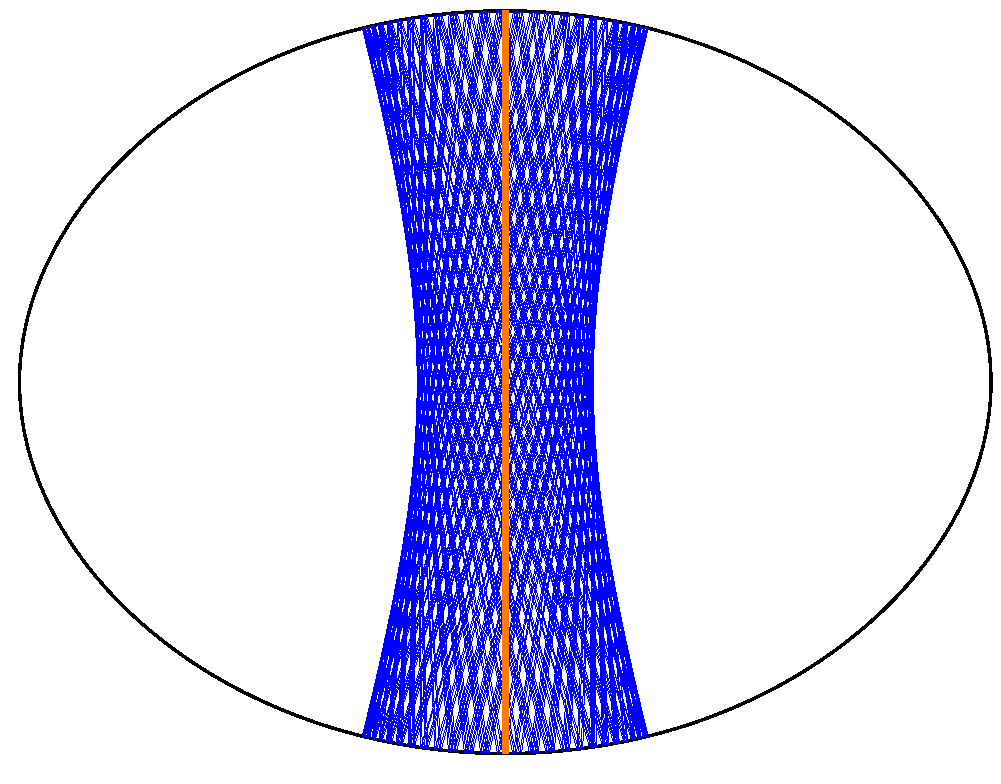}\\ 
\vskip.8cm 
\includegraphics[width=\hsize]{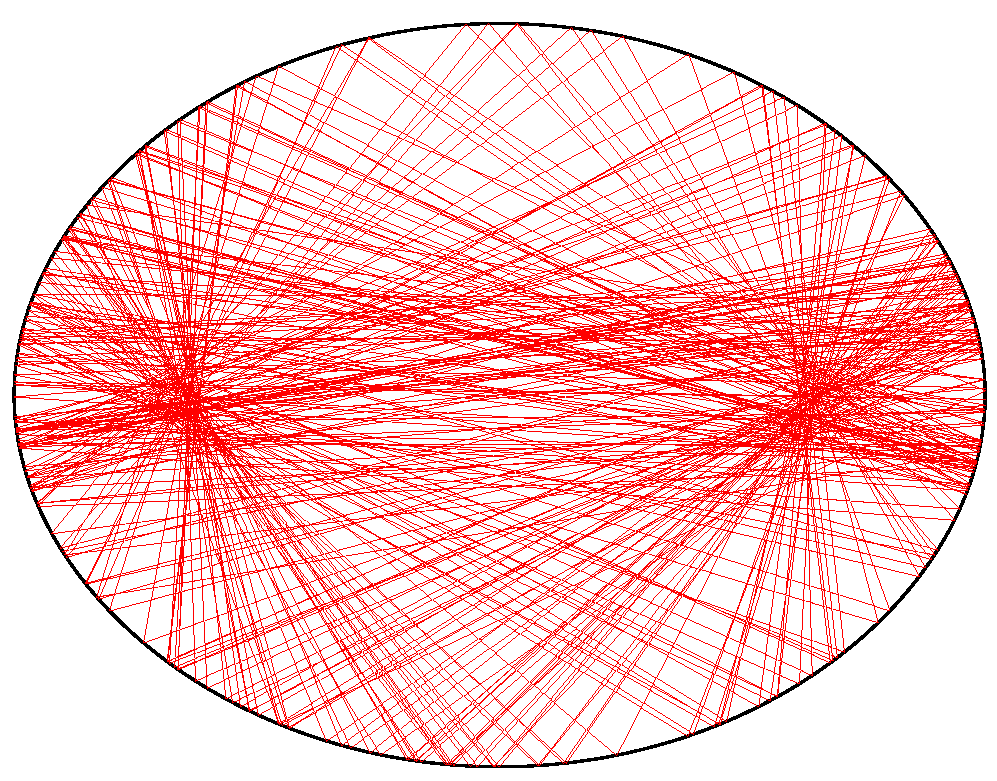}
}
\hskip.5cm 
\parbox{.35\hsize}{
\includegraphics [width=\hsize,height=\hsize]{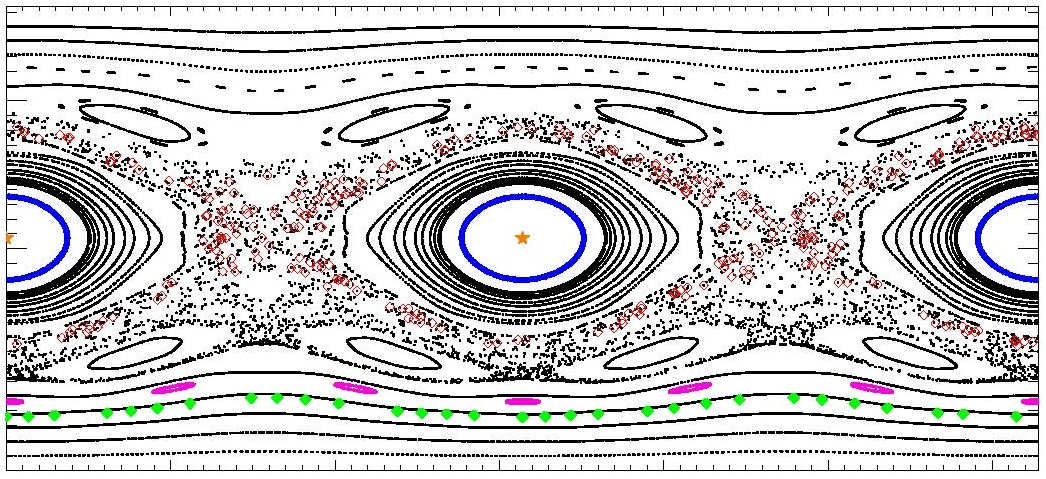}
}
\hskip.5cm\parbox{.14\hsize}{
\includegraphics[width=\hsize]{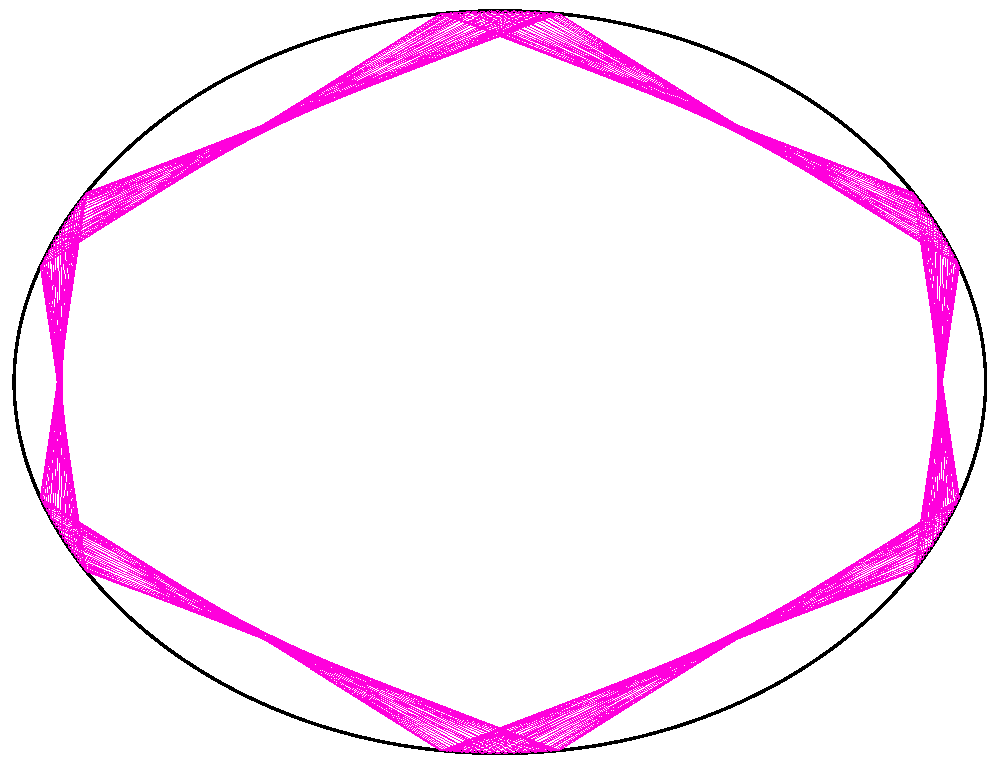}\\ 
\vskip.8cm 
\includegraphics[width=\hsize]{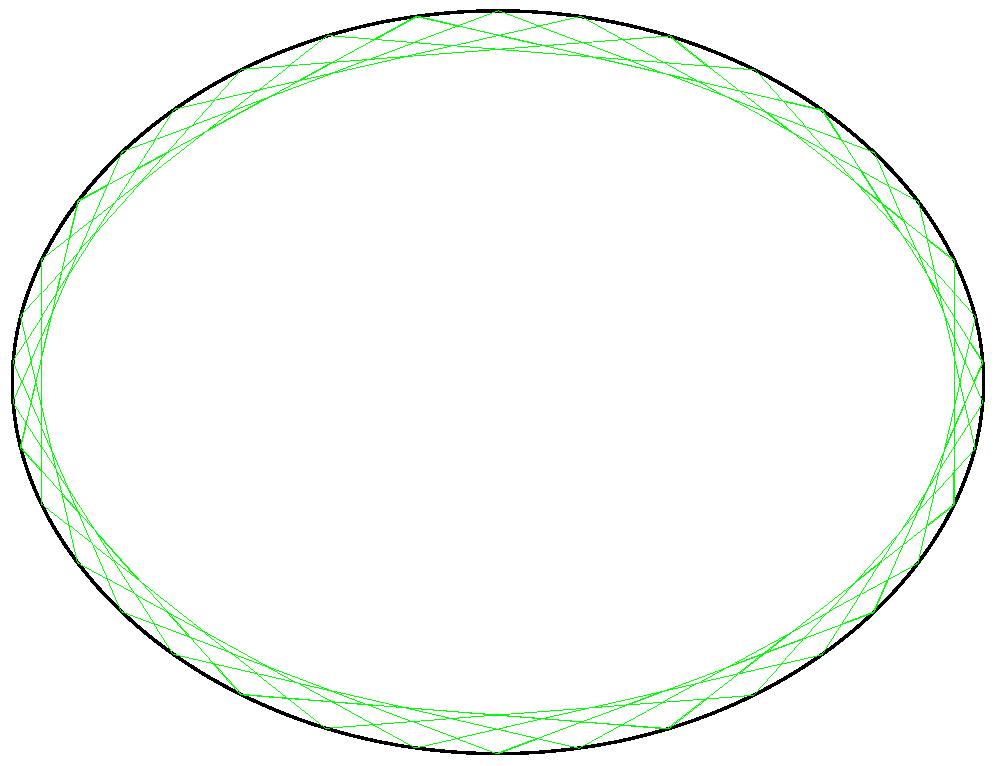}
}
\end{center}
\caption{Numerical simulation of orbits and trajectories of a classic billiard map in a 4-periodic symmetric oval. Blue, pink and green ones correspond to regular motions, inside islands. The red one belongs to the instability set.}
\label{fig:imita}
\end{figure}

These two ingredients give rise to a so called optical billiard formulation, where the light travels inside a closed region along the optical path defined  by the medium refractive index and undergoes elastic reflections at the impacts with the boundary.
Integrating numerically the equations of motion, for a given frequency of rotation, Ligni\`eres and Georget 
 obtained a pictorial  representation of the phase space \cite{lig},
and observe the coexistence of elliptic islands and "chaotic" seas surrounding them.
This coexistence is important for the determination of two subsets of frequency modes: a regular one, corresponding to the islands, and a irregular one, corresponding to the "chaotic" trajectories, that will enable them to describe the pulsation of the star.
The resemblance between these figures and those obtained for classic billiards (see figure~\ref{fig:imita}) naturally leads to the question:
 { \em "Does a billiard dynamics in a refractive medium have the same characteristics as the classic ones?"} 
This is the problem we propose to address here.

We will approach optical billiards as a particular case of billiards defined on general surfaces for which some basic results have already been established \cite{cassio1}.
Thus, in section~\ref{sec:superficies}  we identify a medium 
with a refractive index depending on the distance from the origin 
as a Riemmanian surface with a metric that is a functional multiple of the Euclidean metric. We then describe some geometric properties of this surface and the curves on it.

Billiards are introduced in section~\ref{sec:bilhares} and the description of the dynamics through 
a discrete map of  is presented. The known results for convex billiards in the Euclidean (flat) plane $\mathbb R^2$ (classic billiards) are summarized.
These results allow a generic description of the dynamics and the structure of the phase space.
In the remaining of this section, we focus on optical billiards. 
It was proved in \cite{cassio1} that  geodesically convex billiards on quite general surfaces are well defined by  conservative twist maps acting on a limited cylinder.
The goal here is to specify in which circumstances this is true for optical billiards, the point being the regularity of 
of the surface (determined by the refractive index), the notion of convexity, and in particular the existence and uniqueness of geodesics.
In this case, known dynamical properties of twist maps hold. 

Finally in section~\ref{sec:exemplos} we present some examples of optical billiards. 
We start discussing optical billiards on geodesic circles. It is well known \cite{bialy93,bialy2013hopf} that on surfaces of constant Gaussian curvature, the billiard map is totally integrable  (i.e. the cylindrical phase space is foliated by rotational invariant curves) if and only if it is a (geodesic) circle. 
 For optical billiards, the radial symmetry implies that circles centered at the origin  (which are also geodesic circles) are totally integrable. However, this is probably not true for the other geodesic circles, not  centered at the origin.
 
As a second example we investigate optical billiards on non trivial surfaces of zero Gaussian curvature
(other than the flat plane).
As these surfaces have a singularity at the origin, the results of  \cite{cassio1} do not apply. Nevertheless we can have some insights of the dynamics by the limiting behavior of the geodesics.

At last, we would like to point that {\em refractive billiards} (with refractive instead of reflective boundaries) were also investigated in connection to the study of galaxies \cite{deblasi2023some,barutello-refrac}.
This is a different class of billiards. On the other hand, billiards with a potential acting on the {\em particle} and so affecting its the dynamic inside the region
\cite{debiasi25}  may be connected to optical billiards.

\section{Refracting media as Riemmanian surfaces}

\label{sec:superficies}

Let us consider a two dimensional medium and suppose that, in polar coordinates ${\tt x} = (r \cos \theta, r \sin \theta)$, 
its refractive index is given by a strictly positive function   $n=n(r)$. The optical Lagrangian is then given by 
 $L(r,\theta,\dot r,\dot \theta)=n^2(r)\left({\dot r}^2+r^2{\dot\theta}^2\right)$
and the optical rays will be the minimizers of the associated action, i.e. the solutions of the Euler-Lagrange equations.
\begin{eqnarray}
&& \ddot \theta +  2 \left( \dfrac{1}{r} + \dfrac{n'}{n} \right) \dot r \dot \theta = 0 \label{eqn:E-L}\\
&& \ddot r +  \dfrac{n'}{n}{\dot r}^2 -  \left( \dfrac{1}{r} + \dfrac{n'}{n} \right) r^2 {\dot \theta}^2 = 0 \nonumber
\end{eqnarray}
It is worthwhile to notice that, although the E-L equations are not defined at $r=0$, the half straight lines from the origin ($\dot \theta = 0$) are geodesics. 

In a more geometric context,  the solutions of the Euler-Lagrange equations  are the geodesics of the Riemannian surface $(\mathbb R^2,g)$ defined by the metric  $g(r,\theta)=n^2(r)\left( {dr}^2+r^2{d\theta}^2 \right)$.
As we are using standard polar coordinates, $\{ \dr, \dtheta \}$ 
denote the tangent vectors with respect to the coordinate curves of the parametrization and the coefficients of the First Fundamental Form of the Riemannian surface are 
\begin{equation}
E=    \left< \dr , \dr \right >_g = n^2,   \ \ 
F=  \left< \dr , \dtheta \right >_g = 0,    \ \ 
G= \left < \dtheta , \dtheta \right >_g =n^2 \, r^2
\end{equation}
The metric $g$ and the Euclidenan metric (corresponding to $n=1$) are conformal since 
$\left < u , v \right >_g = n^2  \left < u , v \right >$, which 
implies that angles are preserved, even though lengths are different.  Here $\left< . \,  ,  . \right>$ denotes the usual scalar product of
the flat plane $\mathbb R^2$.

The Christoffel symbols are given by
\begin{equation}\label{eqn:simbolos}
\Gamma_{rr}^r  
= \dfrac{n'}{n}, \
\Gamma_{r\theta}^{\theta} 
= \dfrac{1}{r} + \dfrac{n'}{n}, \
 \Gamma_{\theta\theta}^r  
= - r^2 \left(  \dfrac{1}{r} + \dfrac{n'}{n} \right),  \
 \Gamma_{rr}^{\theta} = 
 \Gamma_{r\theta}^r  = 
 \Gamma_{\theta\theta}^{\theta} 
= 0
\end{equation}
and the  Gaussian curvature of  $(\mathbb R^2,g)$ is 
\begin{equation}\label{eqn:K}
K
= - \dfrac{1}{\sqrt{EG}} \left( \left( \dfrac{E_\theta}{2\sqrt{EG}}\right)_\theta + \left( \dfrac{G_r}{2\sqrt{EG}}\right)_r \right) 
=-\dfrac{1}{n^3} \left( {n''}  + \dfrac{n'}{r} - \dfrac{(n')^2}{n} \right) 
\end{equation}
In particular, if $n$ is at least $C^4$, the curvature $K$ is zero if and only if $n(r) = r^{\beta}$  and, by Minding's Theorem \cite[sec. 4.6]{manfredo} ,   $({\mathbb R}^2, g) $ is then locally isometric to the flat plane. However, for $\beta \ne 0$ the metric has a singularity at the origin.

\subsection{Curves in $({\mathbb R}^2,g)$}
\label{sec:curvas}

\newcommand{\nablag}{{\nabla} }

Let $\gamma (t) = (r(t), \theta(t))$ be a regular $C^2$ parameterized curve in ${\mathbb R}^2$. Its first derivative
$\dot \gamma(t) = \dot r \dr + \dot \theta \dtheta $
does not depends on the metric. On the other hand, the second derivative is given by the covariant derivative
$$
D_g \dot \gamma = \ddot r \dr + \ddot \theta \dtheta + 
\underbrace{ 
\dot r^2 \nablag_r^r +2 \dot r \dot \theta  \nablag_r^{\theta} + {\dot \theta}^2  \nablag_{\theta}^{\theta}
}_{\hbox{\scriptsize connection of the metric $g$}}
$$
where in terms of the Christoffel symbols
\begin{eqnarray}
\nablag_r^r &= \Gamma_{rr}^r \dr + \Gamma_{rr}^{\theta} \dtheta 
&= \dfrac{n'}{n} \dr  \nonumber \\
\nablag_r^{\theta} &= \Gamma_{r\theta}^r \dr + \Gamma_{r\theta}^{\theta} \dtheta
&= \left( \dfrac{1}{r}+  \dfrac{n'}{n} \right) \dtheta \label{eqn:conexao}\\
\nablag_{\theta}^{\theta} &= \Gamma_{\theta \theta}^r \dr + \Gamma_{\theta \theta}^{\theta} \dtheta
& =- r^2 \left(  \dfrac{1}{r} +  \dfrac{n'}{n} \right) \dr \nonumber 
\end{eqnarray}
As the geodesics satisfy $D_g \dot \gamma = 0$, we retrieve the Euler-Lagrange equations \ref{eqn:E-L}.

Let ${\nu} = - r \dot \theta \dr + \dfrac{\dot r}{r} \dtheta $ be a normal vector of the curve, i.e. 
$\left< {\nu}, \dot \gamma \right > = \left< {\nu} , \dot \gamma \right >_g = 0$. 
Moreover $\nu$ and $\dot \gamma$ have the same norm as
\begin{eqnarray*} 
 &\lVert \,   \dot \gamma \rVert_g^2 = n^2 \lVert \,   \dot \gamma \rVert^2 &= n^2 (\dot r^2 + r^2 \dot \theta^2)
\\
 &\lVert  {\nu} \rVert_g^2 = n^2  \lVert  {\nu} \rVert^2 &= n^2 (r^2  \dot \theta^2 +  \dot r^2)
\end{eqnarray*}
(In particular, if we assume that $\gamma$ is parametrized by Euclidean arc length,  we have
$ \lVert \dot \gamma \rVert = \lVert  {\nu} \rVert = 1$. )

The geodesic curvature of $\gamma$ with respect to the metric $g$ is then defined by
$$
\left <{\nu} , D_g \dot \gamma \right >_g =  \kappa_g \, \lVert  {\nu} \rVert_g   \lVert \,  \dot  \gamma \rVert_g^2 
$$


The curvature of $\gamma$ in the {flat plane}
is given by
$
\left <{\nu}, D \gamma' \right > = \kappa \, \lVert  {\nu} \rVert   \lVert \, \dot  \gamma \rVert^2 
$
where the Euclidean connection $D$ is obtained from \ref{eqn:conexao} by taking $n(r) = 1$.
We want to relate $\kappa_g$ to $\kappa$, the point here being that $D_g \ne D$. 
As 
$
\kappa_g =\dfrac{ \left <{\nu} , D_g \dot \gamma \right >_g }
{ \lVert  {\nu} \rVert_g   \lVert \,  \dot  \gamma \rVert_g^2 }
= \dfrac{ n^2 \left <{\nu} , D_g \dot \gamma \right > }
{ n^3  \lVert  {\nu} \rVert   \lVert \,  \dot  \gamma \rVert^2 }
$
we can write
$$
n \kappa_g = \dfrac{ \left <{\nu} , D_g \dot \gamma \right > }
{  \lVert \,  \dot  \gamma \rVert^3 }
=\kappa
+
\dfrac{ \left <{\nu} , (D_g - D)\dot \gamma \right > }
{  \lVert \,  \dot  \gamma \rVert^3 }
$$
Now
$$ (D_g - D)\dot \gamma 
=\dfrac{n'}{n} \left(  \dot r^2 \dr + 2 \dot r \dot \theta  \dtheta - \dot \theta^2 r^2 \dr   \right) 
=  \dfrac{n'}{n} \left(
\dot r \dot \gamma+ r \dot \theta {\nu}  
 \right) 
$$
So 
$$
{ \left <{\nu} , (D_g - D)\dot \gamma \right > }
=   \dfrac{n'}{n} r \dot \theta   \lVert \,  \dot  \gamma \rVert^2 
$$
and finally 
\begin{equation} \label{eqn:kg} 
\kappa_g = \dfrac{\kappa}{n} +  \dfrac{n'}{n^2} \dfrac{r \dot \theta}{  \lVert \,  \dot  \gamma \rVert}
\end{equation}
If $\gamma$ is a closed simple curve enclosing the origin we can use $\theta$ as the parameter and  $\dot \theta = 1$.
This implies that if the refractive index is an increasing function, curves with positive curvature 
with respect to the Euclidean metric will also have positive geodesic curvature with respect to the metric $g$. 

\section{Billiards}
\label{sec:bilhares}
The model proposed by Ligni\`eres and Georget \cite{lig} for the acoustic modes in rapidly
rotating stars leads to the study of the dynamics of optical billiards by identifying the path of the light ray with the trajectory of a particle in the billiard.

The classic billiard problem consists of describing the dynamics of a particle in a bounded region of the plane.
The particle moves freely through the region and undergoes elastic collisions with the boundary.  
This is a two degrees of freedom problem with an infinite potential well defining the boundary of the region. 
As energy is conserved and  choosing a Poincaré section at the boundary, the  billiard problem can be described by a 
two-dimensional map. The two variables correspond to the point of impact at the boundary and the direction of movement when leaving it.  The billiard map associates one impact to the next one and acts therefor on a cylindrical phase space. 

It is well known that the dynamical properties of a classic billiard map strongly depends on the shape of the boundary. 
In particular non convex boundaries are dispersing and lead in general to chaotic/ergodic behavior. 
On the other hand, convex billiards usually present a coexistence of regular and irregular behavior,
although a complete and rigorous  generic description has not yet been established.

Clearly, the billiard problem can be formulated on more general surfaces where the free motion inside the region follows the geodesics. From this point of view, an optical billiard is given by a closed simple curve in the plane with a metric defined by a refractive index  as in the previous section. The dynamical properties of an optical billiard will then depend on both the boundary curve and the refractive index.



\subsection{Classic oval billiards}
\label{sec:classical}

We will focus on convex billiards and describe here their fundamental dynamical properties. 
Let $\gamma$ be an oval , i.e. a closed, regular, simple, oriented counterclockwise, $C^k$ curve, $k\geq 2$, with strictly positive curvature in $\mathbb R^2$.
The billiard map $T$ defines the dynamics by a sequence of impacts. 
Each collision is given by a point on the the boundary and a direction of motion from it.
It is convenient to use the canonical variables: the arc length parameter of the curve $ s \in S^1$ to specify the point on the boundary curve and the tangential momentum $p  = \cos \alpha \in (-1,1)$ to specified the direction of the out going trajectory, 
$\alpha$ being the angle between the velocity and the tangent to the boundary at the impact point.
Therefor, to each oval $\gamma$ is associated a billiard map $T$ from the cylinder $S^1\times(-1,1)$  into itself which to each initial condition associates the next impact and direction $T(s_0, \cos \alpha_0)=(s_1, \cos \alpha_1)$.

The billiard map has some very well known properties (see for instance  \cite{kat}): it is a reversible  $C^{k-1}$-diffeomorphism
with the monotone twist property and preserves the measure $ds\, dp$.  
The twist property means that the image of a vertical line is a graph. In the case of billiards, this can be transladed into the fact that,
as the outgoing direction varies monotonically, so does the next impact parameter, i.e. ${\partial s_1}/{\partial \alpha_0} > 0$.
Moreover, $T$ has no fixed points but Birkhoff \cite{bir27} proved that it has periodic orbits of any period greater or equal to two.

\begin{figure}[ht]
\begin{center}
\includegraphics[width=.25\linewidth]{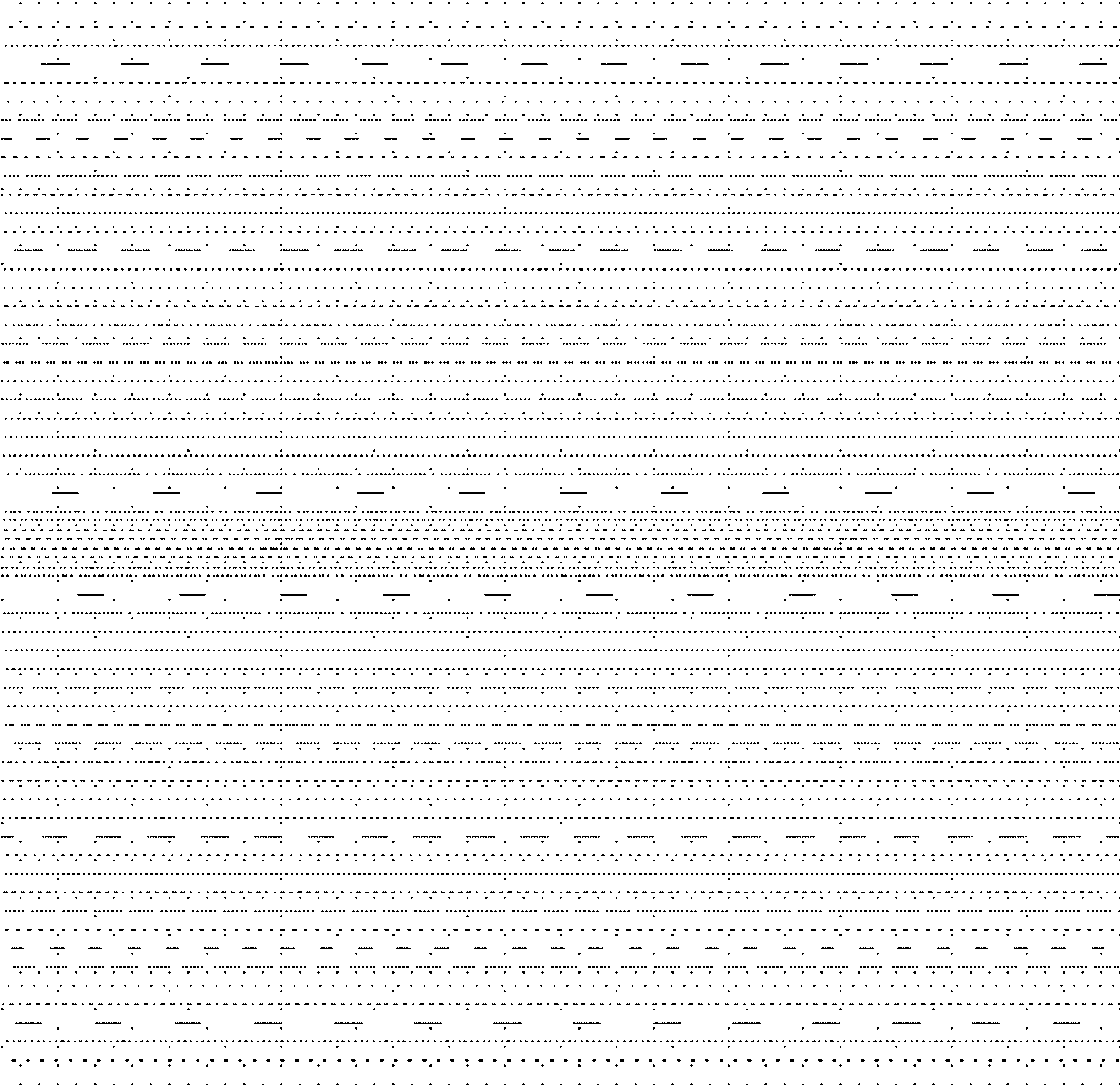}\hskip 2cm
\includegraphics[width=.25\linewidth]{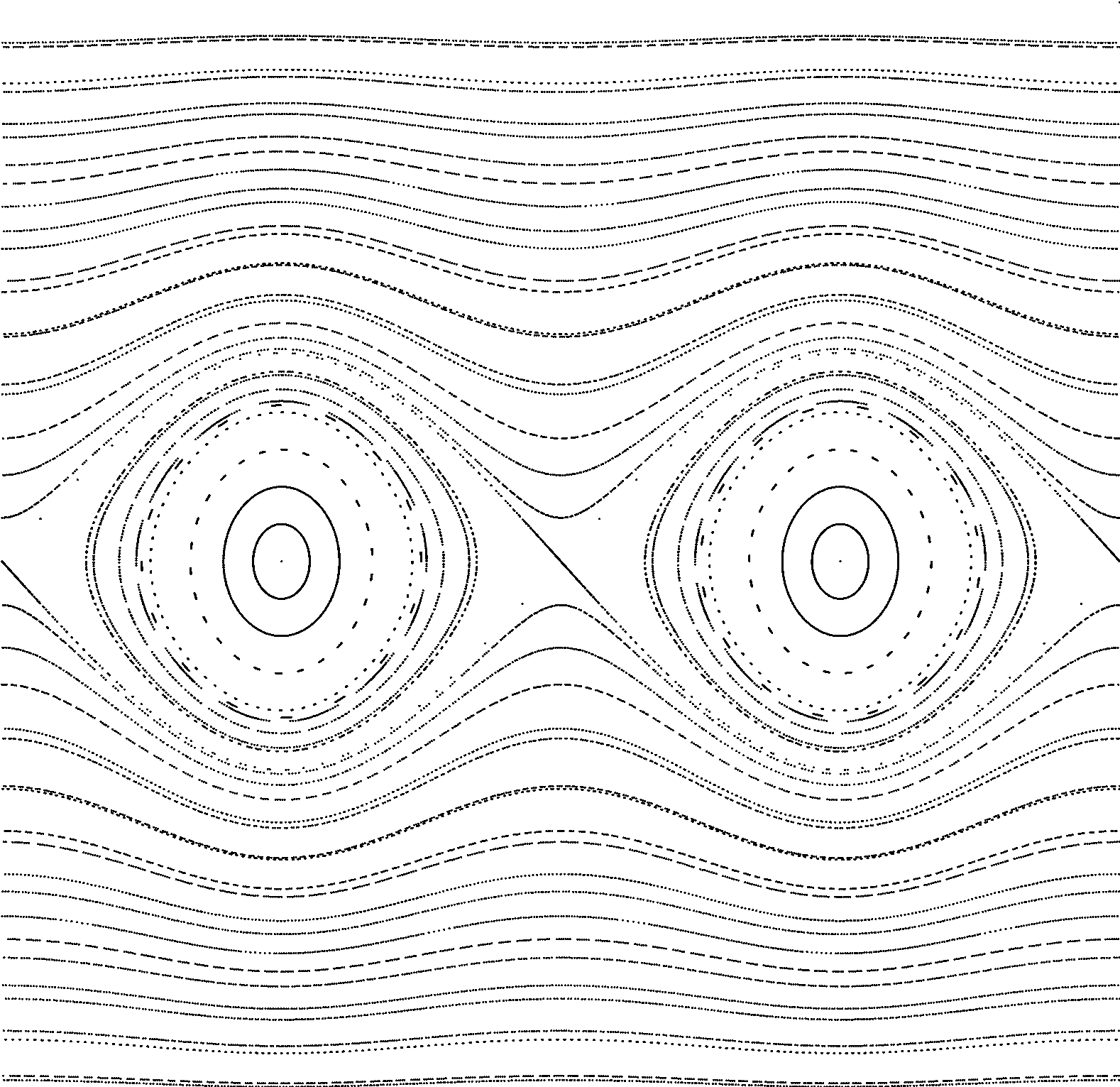}\hskip 2cm
\includegraphics[width=.25\linewidth]{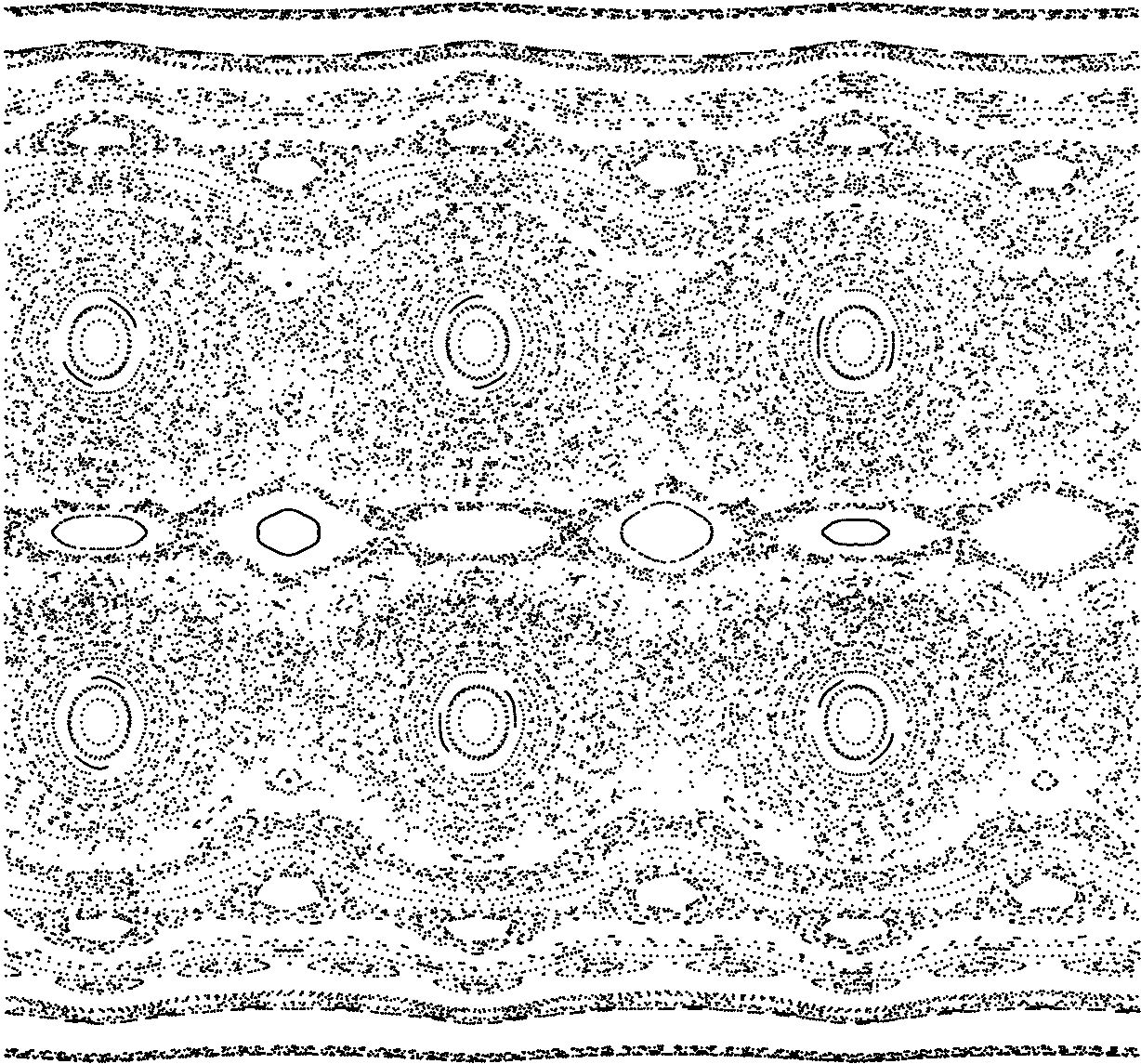}
\end{center}
\caption{Phase space of the billiard in a circle (L), in an elipse (M) and in an oval ().  
}
\label{fig:rot}
\end{figure}
 
A closed, simple, continuous curve in the cylindrical phase space which is not homotopic to a point is called a rotational curve. It is invariant if it is mapped by $T$ onto itself. 
The phase space of the circular billiard, for instance, is foliated by  invariant rotational curves as one can see on figure~\ref{fig:rot} (L).
The circular billiard is said to be totally integrable. 
Elliptical billiards are also integrable (every point of the phase space lies on an invariant curve),
but not totally as some invariant curves are not rotational, but homotopic to an elliptical periodic point (elliptical islands).
A challenging problem about billiards is due to Birkhoff who conjectured that the only integrable convex billiards are ellipses \cite{bir27a,poritsky}. 
Despite many recent advances, this conjecture is still open. 
 
For sufficiently differentiable ovals \cite{laz,dou},  the twist property implies that $T$ has  invariant rotational curves in any small neighborhood of the upper and lower boundaries of the cylinder, as for instance the billiard on figure \ref{fig:rot} (right).
On the other hand, if the boundary has a point of zero curvature, the billiard map has no rotational invariant curves \cite{mather82}. 

Birkhoff \cite{bir32} called the region between two invariant rotational curves, with no other invariant rotational curves inside, an instability region. In \cite{etds} it is shown that for an open and dense set of ovals, the associated billiards always have at least one instability region.
Generically on each instability region one can see the coexistence of 
a countable (maybe finite) number of periodic islands, containing the regular motion, and an instability set, that corresponds to the closure of the stable, or the unstable, curves of the periodic hyperbolic orbits that are outside the islands \cite{dia}. On the instability set, the motion  seems numerically "chaotic", as can be seen on the example pictured bellow.

\subsection{Optical billiards}

The properties mentioned above, which hold for convex billiards on the flat plane, still hold in surfaces with constant Gaussian curvature \cite{luciano,florentin}.
{As for}
more general surfaces, although some basic properties were established in \cite{cassio1}, less is known.
Our goal is to approach optical billiards in this more general context and identify which properties hold in this particular case.
Ultimately, inspired by the figures obtained in \cite{lig} we would like to determine for which refractive indexes the generic dynamic behavior observed in convex Euclidean billiards can be observed.

Together with Dias Carneiro, we proved in \cite{cassio1} that a convex billiard in a totally normal neighborhood of a Riemanniam surface with a sufficiently smooth metric, is a well-defined twist map and is a $C^k$ conservative 
diffeomorphism  if the boundary is  $C^{k+1}$.
In particular, as in the flat case,  if $k\ge 6$ and the geodesic curvature is positive, 
this implies the existence of invariant rotational curves arbitrarily close to the boundary of the phase space. 
Birkhoff's result about the existence of orbits of any period also stands.

The hypothesis about the properties of the region containing the billiard is crucial for the dynamics to be well defined
as the definition of a totally  normal neighborhood implies the existence of a single  geodesic between any two points in it. 
Moreover, convexity implies that this geodesic is entirely  contained in the region.
If the refractive index, is sufficiently differentiable, a totally normal neighborhood of the origin always exists, even though it can be very small. 
However, although conditions to identify totally normal neighborhoods exist \cite{gulliver}, and may be verified in specific cases, it seems difficult estimate its size .

An optical billiard is then defined by given a simple and closed plane curve $\gamma$, and a smooth strictly positive function $n$. We apply the denomination {\em optical billiard} assuming implicitly that the conditions above are met for the region containing the curve $\gamma$ when considering the metric defined by the refractive index $n$.
If the geodesic curvature of $\gamma$  with respect to the metric  is positive, we will say that  the optical billiard is {\em convex}. This implies that the region delimited by the curve is convex, i.e. contained in one of the two half planes determined by a tangent geodesics. 
At this point, it is worth to mention that, although for the original problem of the dynamics of rays in rotating stars, 
optical billiards in convex Euclidean curves should  be considered, these may not be convex when the metrics associated to the refractive is used (equation \ref{eqn:kg}). 

The dynamics in an optical billiard is given, as in the classic case, by a map which associates an impact with the boundary curve $\gamma$ to the next one. Each impact is identified by a point at the boundary and an outgoing direction which, by the assumed normality of the region,  determine a unique geodesic. Convexity implies that this geodesic will intersect the boundary at a single point. The new direction is then defined by the reflexion law.
In what follows, we will give a direct proof that if the boundary is convex and $C^{k+1}$ the optical billiard with a smooth refractive index is a $C^{k}$ diffeomorphism.

\begin{figure}[h]\label{fig:coords}
\begin{center}
\includegraphics[width=.3 \hsize]{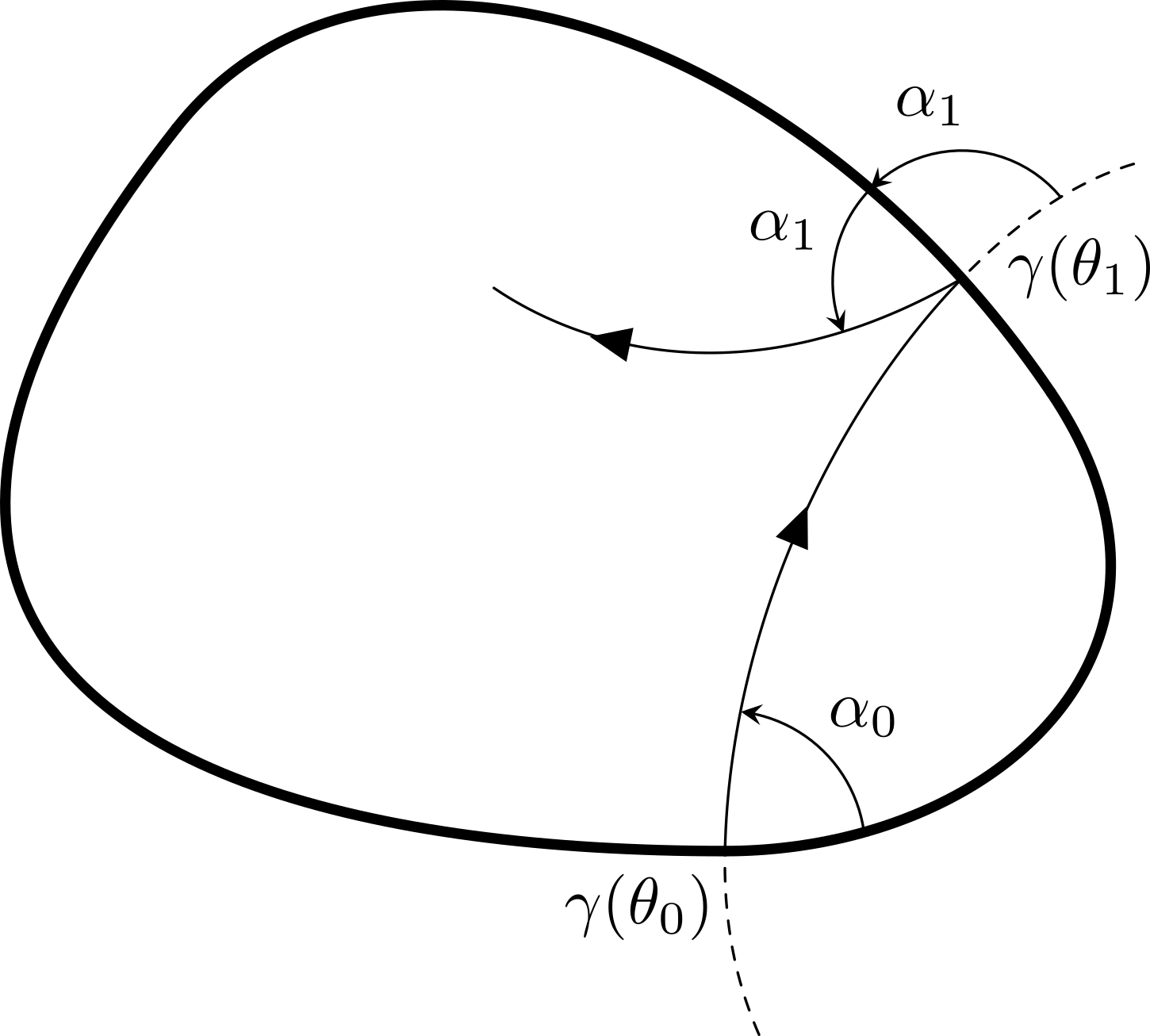}
\end{center}   
\caption{Pictorial description of the optical billiard map.}
\label{definicao}
\end{figure}

Due to the rotational symmetry of $n$, it is useful to use  the standard polar coordinates.  Let the curve be given by $\gamma(\theta) = (\rho(\theta),\theta)$ and, as before,  let $\alpha$ is the angle between the outgoing direction and  the tangent of the curve as on figure~\ref{fig:coords}.
Using these variables, the optical billiard map is denoted by
$$
\begin{array}{ccc}
T:[0,2\pi)\times (0,\pi)&\longrightarrow &[0,2\pi)\times (0,\pi)\\
(\theta_0,\alpha_0)&\longmapsto &(\theta_1,\alpha_1)
\end{array}
$$

The geodesic $\eta(t) = (r(t), \theta(t))$ from $\gamma(\theta_0)$ to $\gamma(\theta_1)$ is the solution of the Euler-Lagrange 
equations (\ref{eqn:E-L}) with initial conditions are given by
$\eta(0) =  \gamma(\theta_0)$ 
and 
$ \angle \left ( \frac{d\eta}{dt} (0), \frac{d\gamma}{d\theta}(\theta_0) \right ) = \alpha_0 $.
If $n$ is smooth then the geodesic $\eta(t;r_0,\theta_0,\dot r_0,\dot\theta_0)$ depends smoothly on $t$ and on the initial conditions $(r_0,\theta_0,\dot r_0,\dot\theta_0)$. 
On the other hand if $\rho$ is $C^{k+1}$ then the initial conditions are $C^k$ on $\theta_0$ and $\alpha_0$. 
Writing the geodesics as $\eta(t;\theta_0,\alpha_0)$ it is smooth on $t$ and $C^k$ on $\theta_0$ and $\alpha_0$.

The parameter $t_1 \ne 0 $ of the next intersection of the geodesic with the boundary $\gamma$ is the solution
of $\rho(\theta(t))=r(t)$. The Implicit Function Theorem implies that $t_1 = t_1(\theta_0,\alpha_0))$ is a $C^k$ function of the initial conditions. Then $\theta_1=\theta(t_1(\theta_0,\alpha_0))$ is also a $C^k$ function of $\theta_0$ and $\alpha_0$.
As the new direction angle $\alpha_1$ is obtained by the reflexion law it follows that if the boundary is $C^{k+1}$, the optical billiard map $T$ is $C^k$. Finally, as the inverse map is obtained by reversing time, $T$ is a $C^k$ diffeomorphism.

As for classic billiards, the invariant measure comes from the conservative character of dynamics and the twist property from the convexity of the curve and the non-existence of conjugate points. We can summarize these results as

\begin{teo}
Consider a medium refractive index given by a smooth strictly positive function $n(r)$ and a  $C^{k+1}$ simple and closed
curve $\gamma$  such that the region limited is contained in a totally normal neighborhood of $({\mathbb R}^2,n) $.
If the geodesic curvature of $\gamma$ is strictly positive, the associated convex optical billiard map is a $C^k$-diffeomorphism with the twist property.
Moreover, if $\gamma$ is given in standard polar coordinates by $r=\rho(\theta)$ the map preserves the measure 
$n(\rho(\theta)) \, \sqrt{\rho'(\theta)^2 + \rho'(\theta)^2} \, \cos \alpha \, d\theta d\alpha$.   
\end{teo}

As it is essentially a consequence of the twist property, the existence of a Lazutkin region (cantor set of invariant rotational curves close to the top and the bottom of the phase space) also follows if $\gamma$ is at least $C^7$ with strictly positive geodesic curvature. 
On the other hand, is follows from \cite{cassio1} that although a point of zero geodesic curvature may not prevent the existence of invariant rotational curves as in the Euclidean case, it certainly prevents the existence of such curves arbitrary close to the top and the bottom of the phase space in optical billiards.

\section{Some examples} \label{sec:exemplos}

\subsection{The geodesic circular billiard}
\label{sec:circulos}

As mentioned in section \ref{sec:classical}, the Euclidean circular billiard is {\it totally integrable}. 
This follows by elementary geometry: as any segment connecting two points of the boundary is the base of an isosceles triangle, the reflexion angle is preserved along  any trajectory and so
the phase space is foliated by invariant rotational horizontal lines. This also can be interpreted as the conservation of the angular momentum.  
On the other hand it is known that the circular billiard is the only totally integrable one \cite{bialy93}. It is also true that a  billiard on  
a constant curvature surface is totally integrable if and only if it is a circle \cite{bialy2013hopf}.

Now, if we consider an optical billiard, the rotational symmetry implies  that circles around the origin are also integrable. 
As before,  we assume that the
circles are contained in a totally normal neighborhood in order to the billiard map be well defined.
We observe that, as  half-lines from the origin are geodesics, the Euclidean circles around the origin are also geodesic circles: 
an Euclidean circle of radius $R$ is a geodesic circle of radius  $\int_0^R n(r)\,dr$.

\begin{teo}\label{thm:circulo}
The optical billiard map in a geodesic circle centered at the origin is totally integrable.
\end{teo}
\begin{proof} 
{ 
Consider a circle of Euclidean radius $R$ centered at the origin and two points on it: $(R,\theta_0)$ and $(R,\theta_1)$. 
Let $\eta_0(t)=(r(t),\theta(t))$ be the geodesic connecting them, i.e. $\eta(0) =  (R,\theta_0)$.  
Since the Euler-Lagrange equations (\ref{eqn:E-L}) do not depend on $\theta$ it follows that 
$\eta_1(t)=(r(t),\theta(t)+\Delta \theta)$, with $\Delta \theta = \theta_1-\theta_0$, is also a geodesic and, by the rotational symmetry, connects 
$(R,\theta_1) = \eta_1(0)$ to $(R,\theta_2) = (R, \theta_1 + \Delta \theta)$
This implies that geodesic triangles with one vertex at the origin and the other two on the (geodesic) circle are isosceles and 
have equal base angles. 
As a consequence, the angle $\alpha$ between the outgoing velocity and the tangent vector is conserved throughout the motion an the billiard is integrable. 
More specifically the horizontal lines of constant $\alpha$ 
in the phase space are invariant.
}
\end{proof}

\begin{figure}[h]
\label{fig:circulo}
\begin{center}
\includegraphics[width=0.25\hsize]{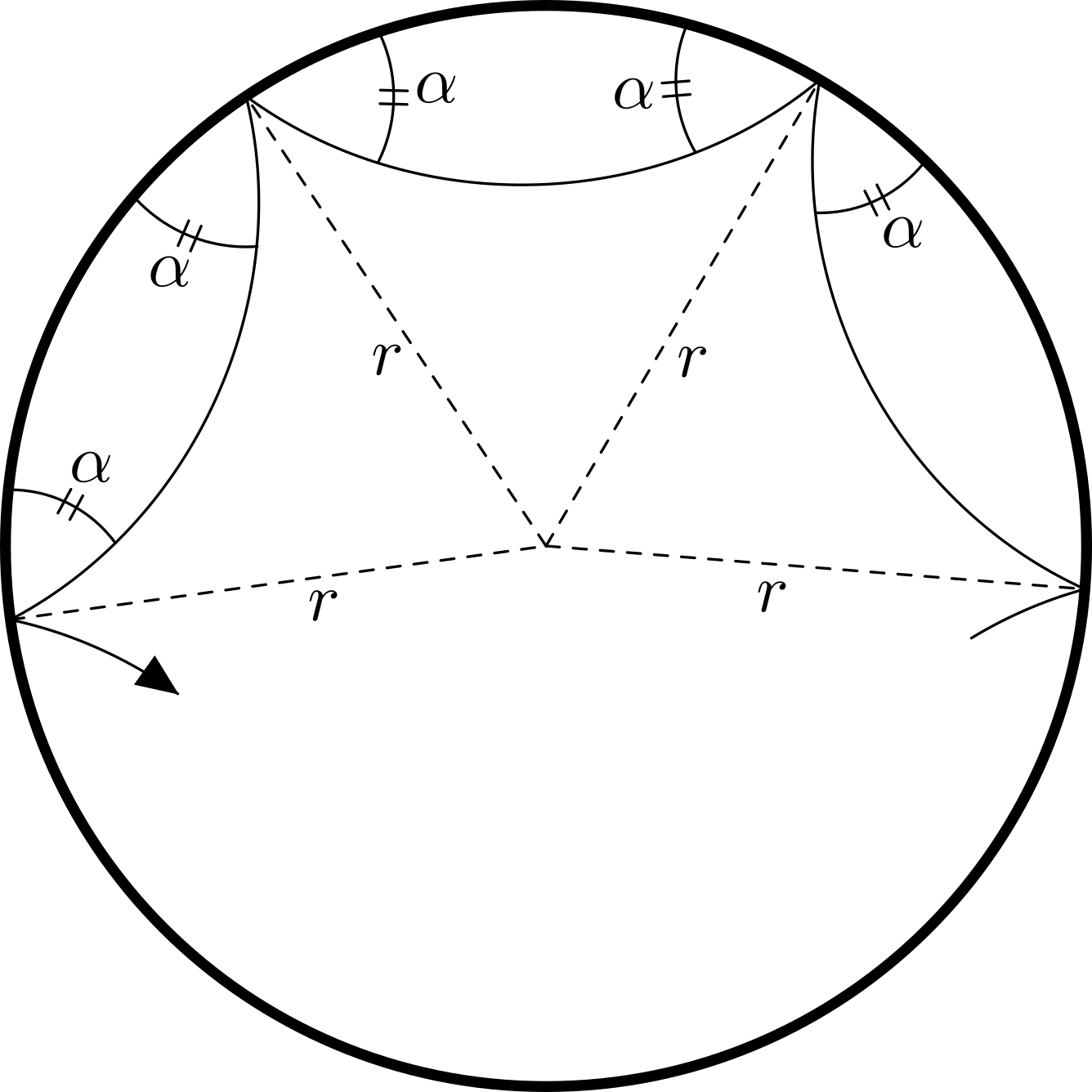}
\end{center}
\caption{Circular billiard}
\end{figure}

A (geodesic) circle of Euclidean radius $R$ has constant curvature $\kappa = \frac{1}{R}$ and if it is centered at the origin, its geodesic curvature
$\kappa_g =  \frac{1}{n(R)}\left(\frac{1}{R} + \frac{n'(R)}{n(R)}\right)$ is also constant. 
However, if $n' <0$, the geodesic curvature of such a circle may be negative. It is worthwhile to notice that the term 
$\left(\frac{1}{r} + \frac{n'(r)}{n(r)} \right)$ also appears in the Christoffel symbols and so its signal plays a significant role 
in the properties of the Riemannian surface  $({\mathbb R}^2,g)$. 
Clearly, a totally normal neighborhood cannot contain a circle of zero curvature, since in this case this circle would itself  be a closed geodesic.  
Thus, this  condition may be used  to estimate the size of the region where the optical billiard is well defined.

As for circles not centered at the origin, Euclidean circles are not geodesic circles and these are not Euclidean circles and do not have constant geodesic curvature. 
Moreover, if we consider a geodesic circle not centered at the origin, 
we do not expect it to be totally integrable as, for general surfaces, the total integrability of a geodesic circle seems to be rather rare
\cite{cassio-circ}.

\subsection{Zero Gaussian curvature refractive media}


It was shown in section \ref{sec:superficies} that if $n(r) = r^{\beta}$ the surface $({\mathbb R}^2, g)$ has zero Gaussian curvature and is therefor locally isometric to the flat plane.
We will discuss the implications of this equivalence on the dynamics.

Let $(r,\theta)$ be the polar coordinates of $({\mathbb R}^2, g) $ and $(R,\Theta)$ the polar coordinates of the flat plane ${\mathbb R}^2$. The respective metrics are $ds^2 = n(r)^2 (dr^2 + r^2 d\theta^2)$ and
$dS^2 = (dR^2 + R^2 d\Theta^2)$, which as we noted before are conformal for any index $n(r)$.
Assuming that $n(r) = r^{\beta}$  it is easy to verify that, if $\beta \ne -1$, the map from 
$({\mathbb R}^2, g) $ to ${\mathbb R}^2$ defined by
$$ R = \frac{r^{\beta+1}}{\beta+1}  \qquad   \Theta = (\beta+1) \theta$$
is an isometry since
$ dS^2 =  r^{2\beta} \, ds^2  $.
This isometry is clearly only local, as it is not a bijection (it is in fact one to $|\beta+1|$). Moreover, as we already mentioned,  the metric is singular at the origin. 

The billiard defined in any closed curve contained in a preimage of 
${\mathbb R}^2-\{0\}$
will thus exhibit the exactly same dynamics of the corresponding curve in the flat plane.
In particular geodesic circles are totally integrable, geodesic ellipses are integrable and the generic properties of geodesically convex curves hold.

As for curves enclosing the origin, the isometry does not apply and so the dynamical equivalence does not hold as it. 
In what follows we will analyze the case $\beta = 1$ which serves as a prototype to $\beta \ge 1$ and show that there is still a correspondence between the dynamics on the flat plane and on the refractive medium.

When $\beta = 1$, the isometry is written in cartesian coordinates as
\begin{eqnarray*}
&&X = R \cos \Theta = \dfrac{r^2}{2} \cos 2 \theta =  \dfrac{r^2}{2} \left (      \cos^2 \theta - \sin^2 \theta \right) =
\frac{x^2}{2} -  \frac{y^2}{2}
\\
&&Y = R \sin \Theta = \dfrac{r^2}{2} \sin 2 \theta =  \dfrac{r^2}{2} \left (2 \, \sin\theta \, \cos \theta \right) =
xy
\end{eqnarray*}
In particular, the image of a straight line $ a X + b Y = c$ is  $ a(x^2 - y^2)/2 + b xy = c$, which implies that the 
geodesics of $({\mathbb R}^2, g) $ are hyperboles.
We also observe that a half plane in $({\mathbb R}^2, g) $ is isometric to the entire plane ${\mathbb R}^2 $

\begin{figure}[h] 
\begin{center}
\includegraphics[width=0.175\hsize]{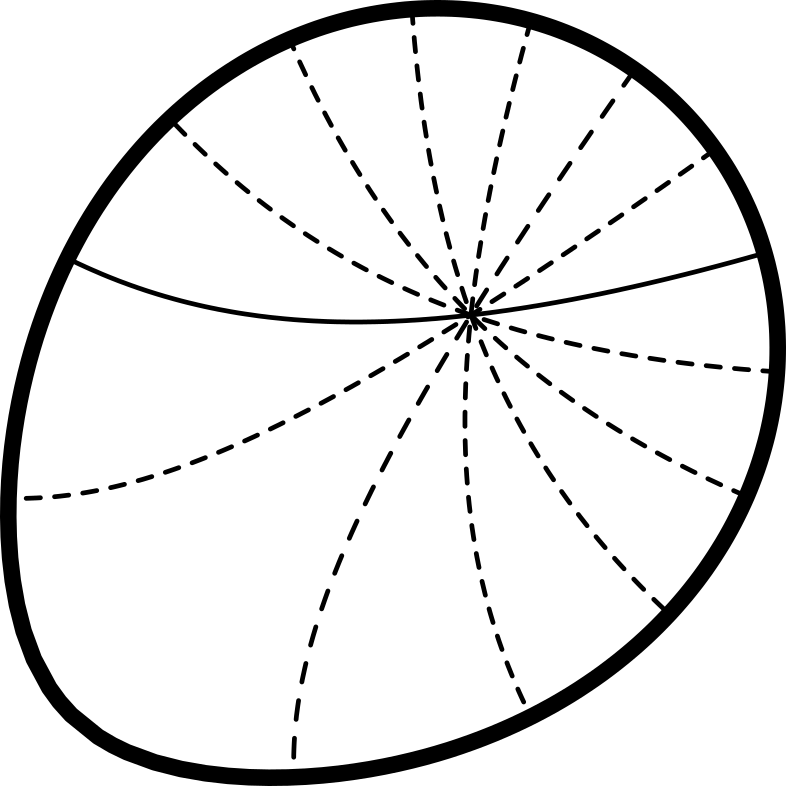} \hskip 1cm
 \includegraphics[width=0.175\hsize]{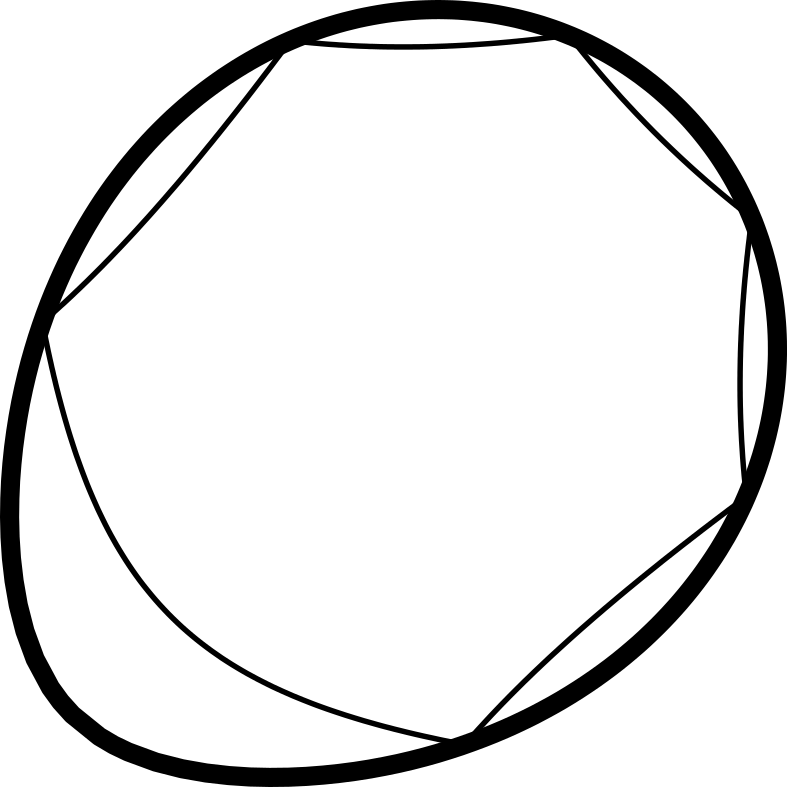} \hskip 1cm
 \includegraphics[width=0.175\hsize]{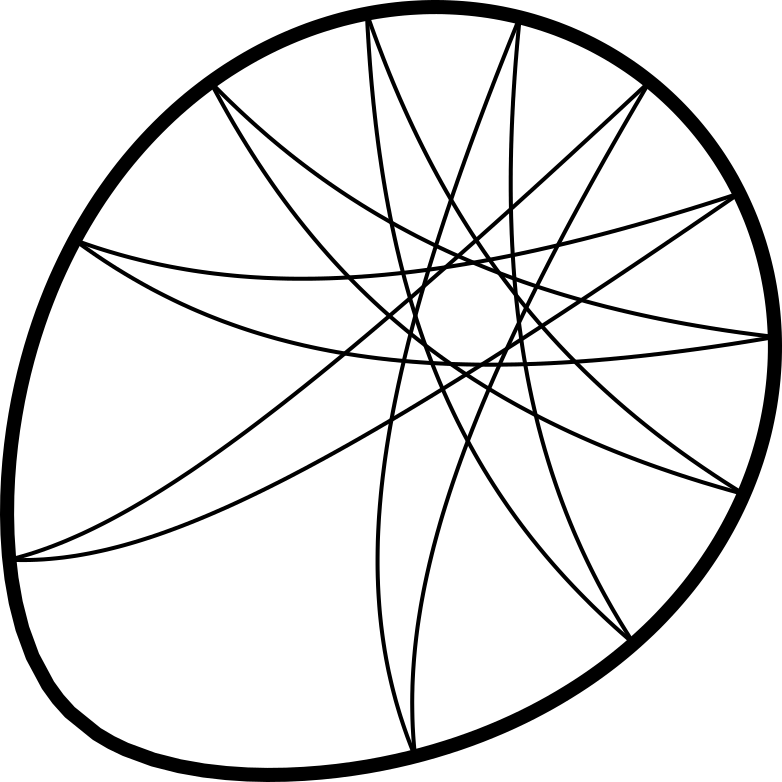} \hskip 1cm 
 \includegraphics[width=0.175\hsize]{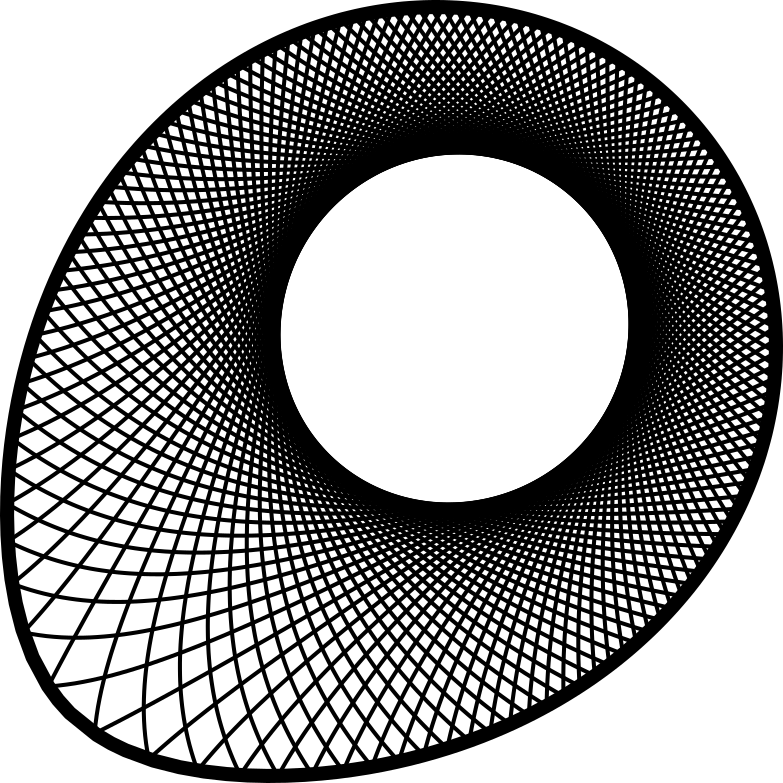} 
\end{center} 
\caption{($n=r$) A geodesic circle off the origin and some trajectories (from left to right: some period 2, a period 6, a period 11 and a dense trajectories).}
\label{fig:r-geocir}
\end{figure}

On figure~\ref{fig:r-geocir} bellow we display a geodesic circle off the origin and some geodesic. We also illustrate some billiard trajectories in it, which, in this case, are mapped to trajectories in the equivalent Euclidean circle.

\begin{figure}[h]
\begin{center}
\includegraphics[width=0.2\hsize]{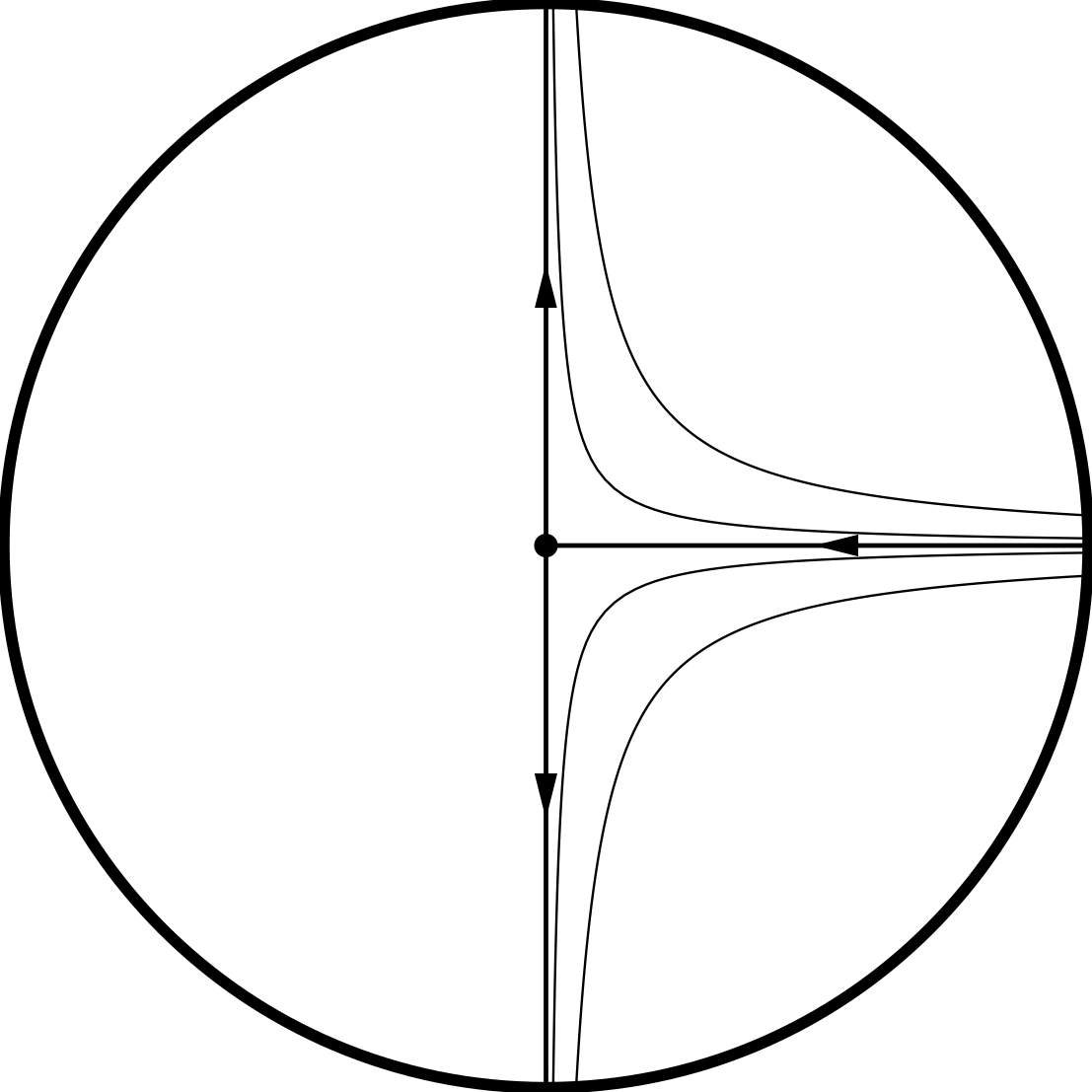} 
\end{center} 
\caption{($n=r$) Non existence of period 2 orbits for the circle around the origin:  discontinous behavior of the geodesics and of the normal trajectories.}
\label{fig:r-singgeo}
\end{figure}

The singularity of the metric at the origin implies a discontinuity of the (straight) geodesics through it as can be observed on 
figure~\ref{fig:r-singgeo}. This implies that a circular billiard centered at the origin in the refractive case does not have period 2 orbits. 
The figure suggest that, as the outgoing direction  converges to the normal direction, the trajectory converges to a (singular) period 4 trajectory.   

In fact, each trajectory in the Euclidean circle of period $p \ge 4$ corresponds by the isometry to a trajectory of period $2p$  in the equivalent circle in 
$({\mathbb R}^2, g) $. This is of course a consequence of the 2-folding of the isometry.

\begin{figure}[h]
\begin{center}
\includegraphics[width=0.175\hsize]{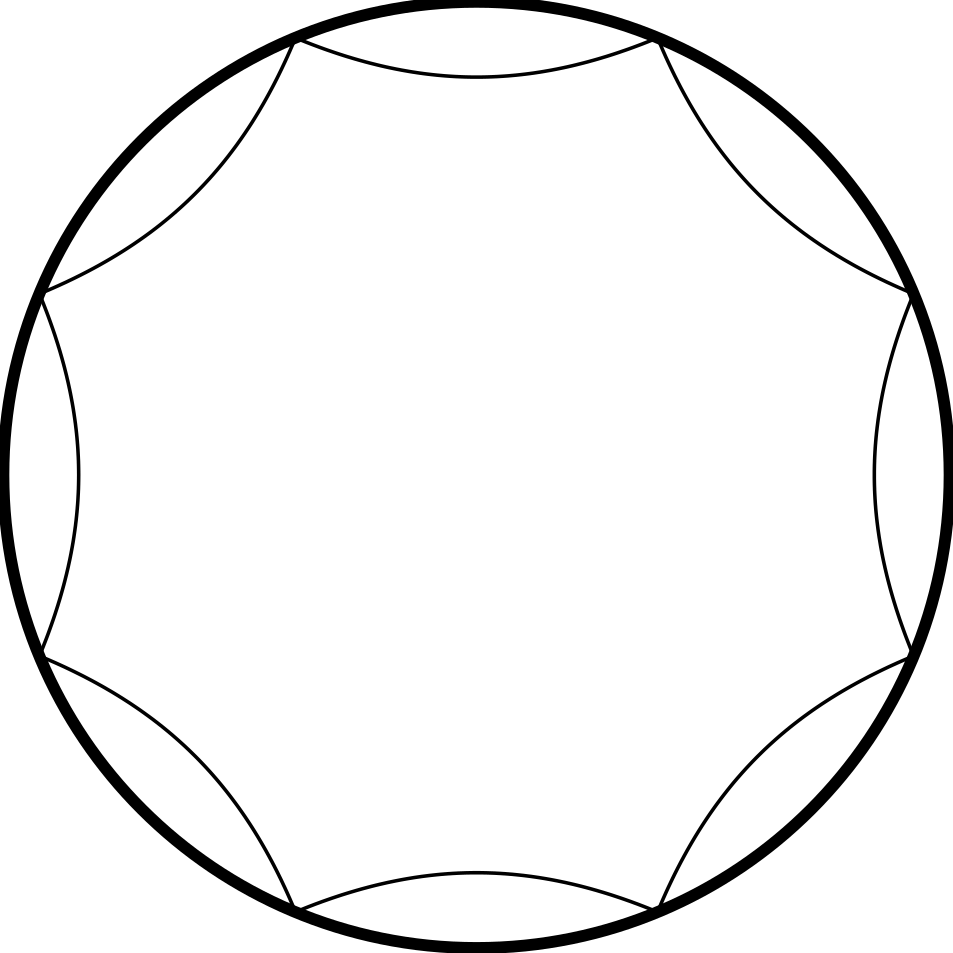} \hskip 1cm
\includegraphics[width=0.175\hsize]{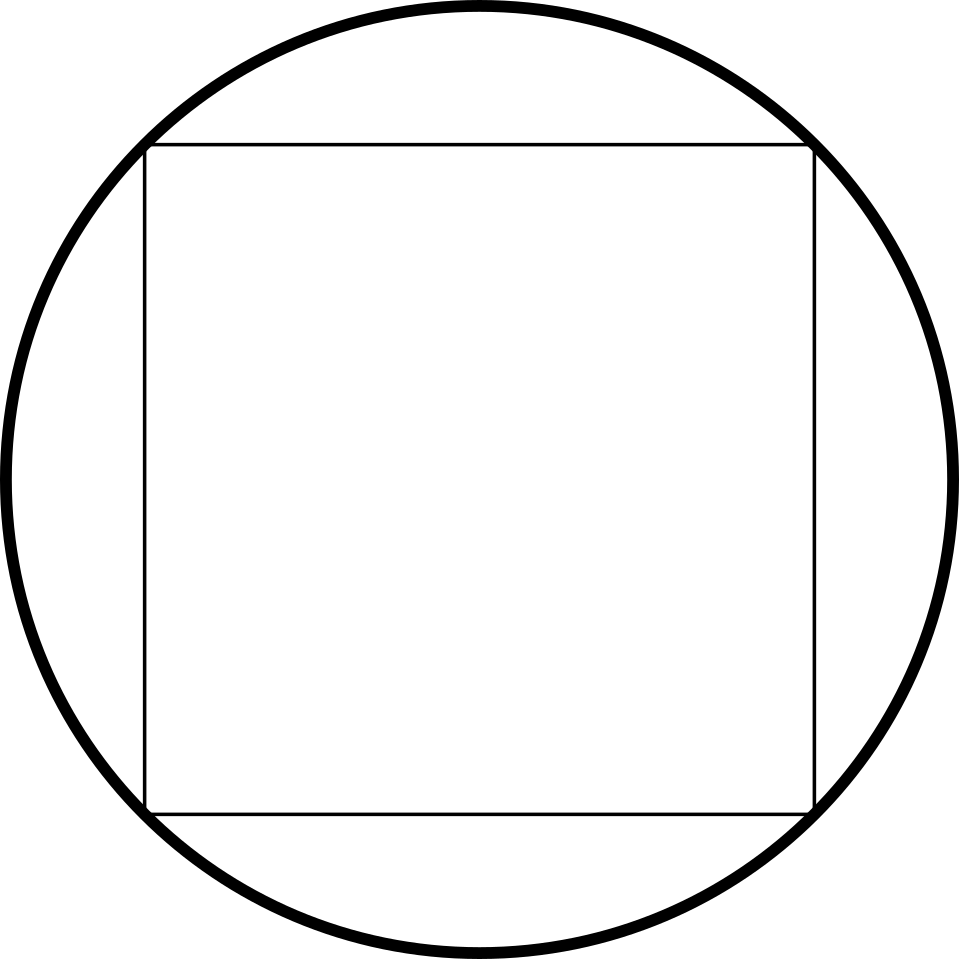}
\end{center} 
\caption{($n=r$) A period 8 trajectory in the refractive case (L) and the corresponding 4 periodic trajectory in the Euclidean case, traveled twice (R).}
\end{figure}

Another interesting example is  $\beta = -2$, i.e. $n(r) = 1/r^2$. In this case the isometry is given by the inversion
$$R = \frac{1}{r} \qquad \Theta = \theta$$ and so the geodesics are circles through the origin which is a highly singular point, as corresponds to the infinite by the isometry.

\begin{figure}[h]
\begin{center}
\includegraphics[width=0.175\hsize]{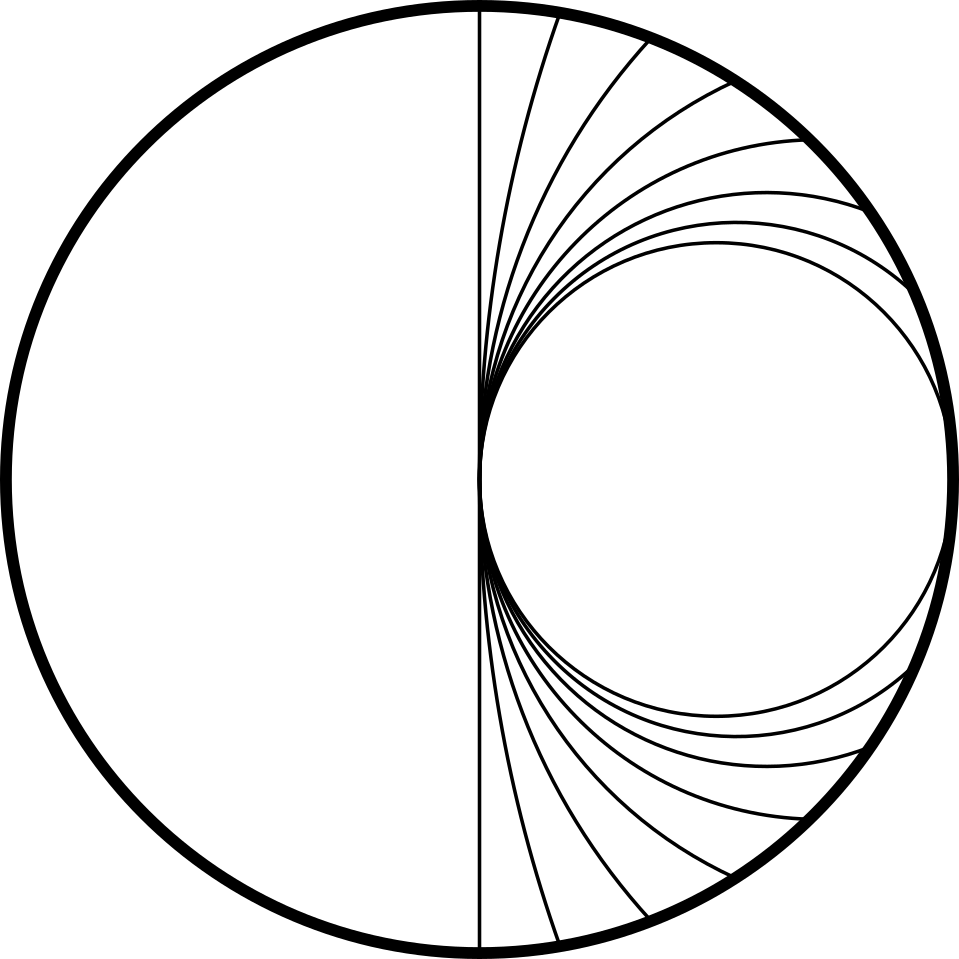} \hskip 1cm
\includegraphics[width=0.175\hsize]{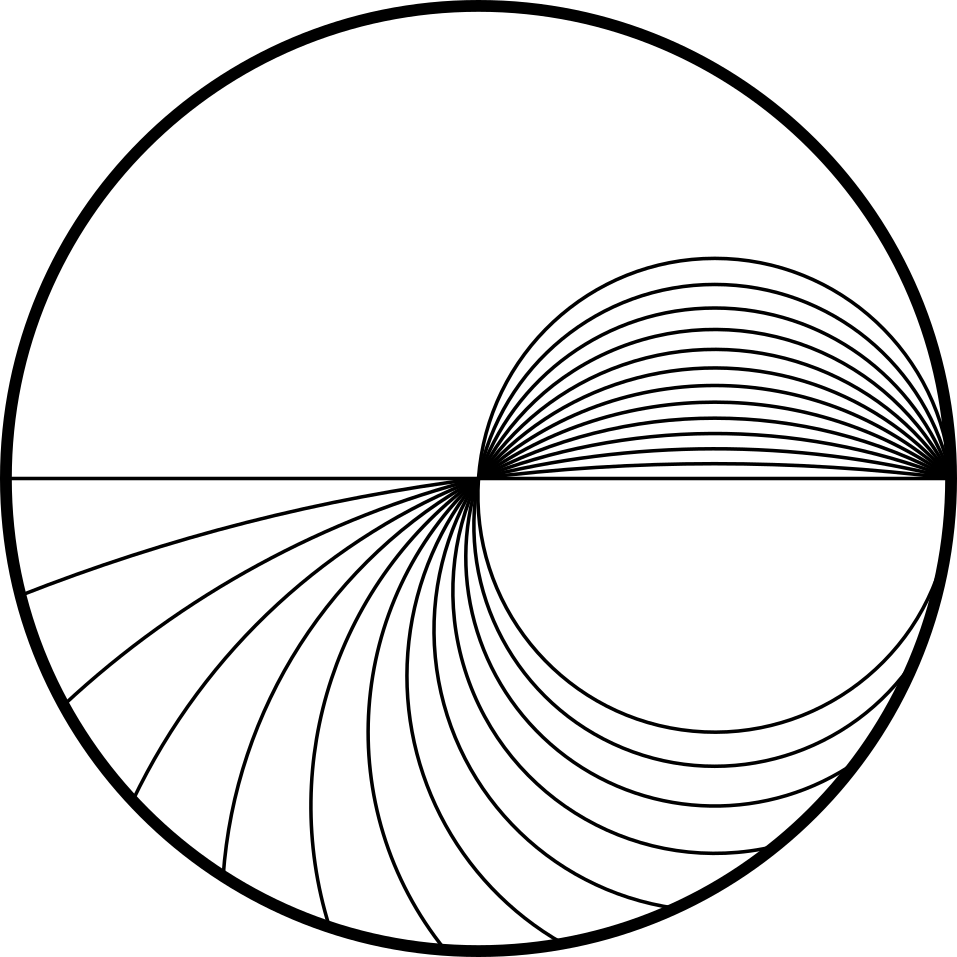} \hskip 1cm
\includegraphics[width=0.175\hsize]{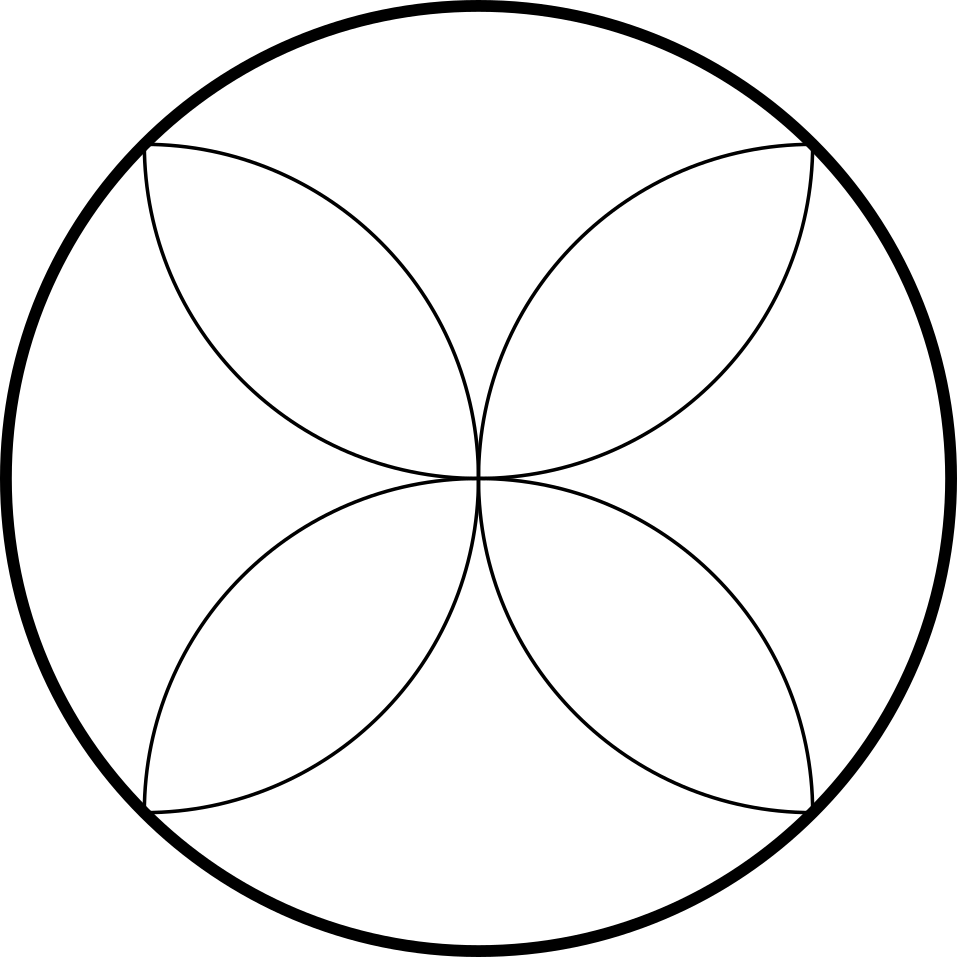} 
\end{center} 
\caption{($n=1/r^2$) Geodesics corresponding to vertical lines (L) , different outgoing directions (M) and a period 4 trajectory (R) for a circle around the origin.}
\end{figure}

\begin{ack}
We thank Daniel R. Reese (Obs. Paris) for drawing our attention to the similarity between the images of rotating stars and those of convex billiards.
\end{ack}

\bibliographystyle{abbrv}

\end{document}